\documentclass[12pt]{article}
\usepackage{amsmath,amssymb,amsthm,amscd}

\numberwithin{equation}{section}

\newcommand{\Z}{{\mathbb Z}}

\newcommand{\C}{{\mathbb C}}
\newcommand{\Sring}{A}
\newcommand{\sK}{K_0}
\newcommand{\sR}{R_0}
\newcommand{\sa}{a_0}
\newcommand{\bo}{b_0}

\DeclareMathOperator{\End}{End} 
\DeclareMathOperator{\Gal}{Gal} 
 
\DeclareMathOperator{\Mat}{Mat} 
 
\DeclareMathOperator{\id}{id}

\DeclareMathOperator{\Top}{top} 
\DeclareMathOperator{\lc}{lc} 
\DeclareMathOperator{\ld}{ld}

\newtheorem{lemma}{Lemma}
\newtheorem{theorem}{Theorem}
\newtheorem{proposition}[lemma]{Proposition}

\theoremstyle{definition}
\newtheorem{definition}{Definition}

\title{Finite-dimensional vertex algebra modules over fixed point commutative subalgebras}
\author{Kenichiro Tanabe\footnote{Partially supported by JSPS
Grant-in-Aid for Scientific Research No. 20740002.}\\\\
Department of Mathematics\\
Hokkaido University\\
Kita 10, Nishi 8, Kita-Ku, Sapporo, Hokkaido, 060-0810\\
Japan\\\\
ktanabe@math.sci.hokudai.ac.jp}
\date{}
\begin{document}
\maketitle
\begin{abstract}
Let $A$ be a connected commutative $\C$-algebra with derivation $D$,
$G$ a finite linear automorphism group of $A$ which preserves $D$,
and $R=A^G$ the fixed point subalgebra of $A$ under the action of $G$.
We show that 
if $A$ is generated by a single element as an $R$-algebra and
is a Galois extension over $R$ in the sense of 
M. Auslander and O. Goldman, then
every finite-dimensional vertex algebra 
$R$-module has a structure of twisted vertex algebra $A$-module.
\end{abstract}
\section{Introduction}
Vertex algebras
and modules over a vertex algebra
were introduced by Borcherds in \cite{B}.
As an example,
every commutative ring $A$ with an arbitrary derivation $D$
has a structure of vertex algebra,
and every ring $A$-module naturally becomes
a vertex algebra $A$-module. 
However, this does not imply that
ring $A$-modules and vertex algebra $A$-modules
are same.
In fact, a vertex algebra $\Z[z,z^{-1}]$-module 
which is not a ring $\Z[z,z^{-1}]$-module was given in \cite[Section 8]{B},
where $\Z[z,z^{-1}]$ is the ring of Laurent polynomials over $\Z$.
This tells us that in general these two kind of $A$-modules are certainly different.
{From} now on, for a commutative $\C$-algebra $A$ with derivation $D$,
we shall call a vertex algebra $A$-module
a {\em vertex algebra $(A,D)$-module} to distinguish 
it from ring $A$-modules.
It is a natural first step to investigate vertex algebra $(A,D)$-modules
to understand modules over general vertex algebras.
In \cite{T1,T2}
for the polynomial ring $\C[s]$ and the field of rational functions $\C(s)$, 
the finite-dimensional vertex algebra modules
which are not $\C$-algebra modules are classified.

Let $A$ be a commutative $\C$-algebra with derivation $D$, 
$G$ a finite linear automorphism group of $A$ which preserves $D$,
and $R=A^G$ the fixed point subalgebra of $A$ under the action of $G$.
In this paper, we shall investigate a relation between vertex algebra $(R,D)$-modules 
and twisted vertex algebra $(A,D)$-modules.
In Theorem \ref{theorem:correspondence},
I shall show that
if $A$ is a connected commutative $\C$-algebra generated by a single element as an $R$-algebra and
is a finite Galois extension over $R$ in the sense of \cite[p.396]{AG},
then 
every finite-dimensional indecomposable vertex algebra $(R,D)$-module becomes 
a $g$-twisted vertex algebra $(A,D)$-module for some $g\in G$.
This is a generalization of \cite[Theorem 1]{T2} and
is related the following open conjecture on vertex operator algebras:
let $V$ be a vertex operator algebra and $H$ a finite automorphism group of $V$.
It is conjectured that under some conditions on $V$,
every irreducible module over the fixed point vertex operator subalgebra $V^H$
is contained in some irreducible $h$-twisted 
$V$-module for some $h\in H$ (cf.\cite{DVVV}).
The conjecture is confirmed for some examples in \cite{AD,DLTYY,DN1,DN2,DN3,TY1,TY2}.
However $A$ is not a vertex operator algebra except in 
the case that $D=0$ and $\dim_{\C}A<\infty$,
Theorem \ref{theorem:correspondence} implies that
the conjecture holds for 
all finite-dimensional vertex algebra $R$-modules in a stronger sense.

This paper is organized as follows:
In Section 2 we recall 
some notation and properties 
of Galois extensions of rings, vertex algebras and their modules.
In Section 3 we show that 
every finite-dimensional indecomposable vertex algebra $R$-module 
becomes a $g$-twisted vertex algebra $(A,D)$-module 
for some $g\in G$.
In Section 4 we give the classification of the
finite-dimensional vertex algebra $\C[s,s^{-1}]$-modules
which are not $\C$-algebra $\C[s,s^{-1}]$-modules.
In Section 5 for the $\C$-algebra $A=\C[s,s^{-1}][t]/(t^n-s)$,
which is a Galois extension over $\C[s,s^{-1}]$ with Galois group
the cyclic group of order $n$,
and for all finite-dimensional indecomposable vertex algebra $\C[s,s^{-1}]$-modules $(M,Y_M)$ obtained in Section 4,
we study twisted vertex algebra $(A,D)$-module structures over $(M,Y_M)$.

\section{Preliminary}

We assume that the reader is familiar with the basic knowledge on
vertex algebras as presented in \cite{B,DLM1,LL}.

Throughout this paper all rings and algebras are commutative and associative 
and have identity elements,
$R$ denotes a ring,
$R[Z]$ denotes the polynomial ring in one variable $Z$ over $R$ ,
$G$ denotes a finite group,
$\zeta_p$ denotes a primitive $p$-th root of unity for a positive integer $p$,
and $(V,Y,{\mathbf 1})$ denotes  a vertex algebra.
Recall that $V$ is the underlying vector space, 
$Y(\cdot,x)$ is the linear map from $V$ to $(\End V)[[x,x^{-1}]]$,
and ${\mathbf 1}$ is the vacuum vector.
Let ${\mathcal D}$ be the endomorphism of $V$
defined by ${\mathcal D}v=v_{-2}{\mathbf 1}$ for $v\in V$.
 
First, we recall some results in \cite{AG,CHR,DI,J} for separable algebras
over a ring. 
A ring $R$ is called {\em connected} if $R$ has no idempotent other than $0$ and $1$.
An $R$-algebra $\Sring$ is called {\em separable} if 
$\Sring$ is a projective $\Sring\otimes_{R}\Sring$-module.
An $R$-algebra $\Sring$ is called {\em strongly separable} if it is finitely 
generated, projective, and separable over $R$.
Let us recall the Galois extension of $R$ introduced in 
\cite[p.396]{AG}. 
The following definition, which is equivalent to that in \cite[p.396]{AG},
is given in \cite[Theorem 1.3]{CHR}.

\begin{definition}\label{definition:galois}
Let $\Sring$ be a ring extension of  $R$ and let $G$ be 
a finite group of automorphisms of $\Sring$.
We denote by $A^G$ the fixed point subring of $A$ under the action of $G$.
The ring $\Sring$ is called a {\em Galois extension} of $R$
with Galois group $G$, if the following three conditions hold:
\begin{enumerate}
\item
$\Sring^{G}=R$.
\item 
For each non-zero idempotent $e\in \Sring$ and each $g\neq h$ in $G$,
there is an element $x\in \Sring$ with $g(x)e\neq h(x)e$.
\item
$\Sring$ is a separable $R$-algebra.
\end{enumerate}
\end{definition}
Note that if $A$ is connected, then the condition (2) in Definition \ref{definition:galois} is always satisfied.
It follows from \cite[Theorem 1.3]{CHR} that if $\Sring$ is a Galois extension of $R$,
then $\Sring$ is a strongly separable $R$-algebra.

In \cite[p.467]{J}, A polynomial $P(Z)\in R[Z]$ is called {\em separable} in case $P(Z)$ is monic and 
the factor ring $R[Z]/(P(Z))$ is a separable $R$-algebra.
In this case, $R[Z]/(P(Z))$ is strongly separable
since $R[Z]/(P(Z))$ is a free $R$-module of rank $\deg P(Z)$.
For an $R$-algebra $\Sring$, an element $\theta\in \Sring$ is called a
{\em primitive element} if $\Sring=R[\theta]$, namely $A$ is generated by a single element $\theta$ as an $R$-algebra. 
It is shown in \cite[Theorem 2.9]{J} that
if $A$ is a strongly separable $R$-algebra and if $A$ has a primitive element, then 
there is a separable polynomial $P(Z)$ such that
$A\cong R[Z]/(P(Z))$ as $R$-algebras.

Let $R$ be a connected ring,
$P(Z)\in R[Z]$ a separable polynomial,
and suppose that the factor ring $A=R[Z]/(P(Z))$ is connected and is a Galois extension of $R$ with 
Galois group $G$. Set $\theta=Z+P(Z)\in A$.
Since $A=R[\theta]$, we have $g(\theta)\neq \theta$ for all $g\in G$ without the identity element.
By \cite[Lemma 4.1]{CHR}
and
\cite[Lemma 2.1]{J}, the order of $G$ is equal to $\deg P(Z)$. 
Thus, $G$ acts regularly on
the set of all roots of the polynomial $P(Z)$ in $A$ and hence
$P(Z)=\prod_{g\in G}(Z-g(\theta))$.  
For an $R$-linear homomorphism $f$ from $A$ to an $R$-algebra $B$,
\cite[Lemma 2.1]{J} says that $f(g(\theta))\neq f(h(\theta))$ for all $g\neq h$ in $G$.
This tells us that if $B$ is an integral domain, 
$f$ induces a bijection 
from $\{g(\theta)\ |\ g\in G\}$ to the set of all roots of $f(P(Z))\in B[Z]$ in $B$.
In particular, $f(P(Z))$ has no multiple root.

Next, we recall some results in \cite{B} for a vertex algebra constructed from a 
commutative $\C$-algebra with a derivation.

\begin{proposition}{\rm\cite{B}}\label{proposition:comm-alg}
The following hold:
\begin{enumerate}
\item
Let $A$ be a commutative $\C$-algebra with identity element $1$ and 
$D$ a derivation of $A$.
For $a\in A$, define $Y(a,x)\in(\End A)[[x]]$ by
\begin{align*}
Y(a,x)b&=\sum_{i=0}^{\infty}\dfrac{1}{i!}(D^{i}a)bx^{i}
\end{align*}
for $b\in A$. Then, $(A,Y,1)$ is a vertex algebra.
\item
Let $(V,Y,{\mathbf 1})$ be a vertex algebra such that
$Y(u,x)\in(\End V)[[x]]$ for all $u\in V$.
Define a multiplication on $V$ by $u v=u_{-1}v$ for $u,v\in V$. 
Then, $V$ is a commutative $\C$-algebra with identity element ${\mathbf 1}$ and 
${\mathcal D}$ is a derivation of $V$.
\end{enumerate}
\end{proposition}
Throughout the rest of this section,
$A$ is a commutative $\C$-algebra with identity element $1$  
and $D$ a derivation of $A$.
Let $(A,Y,1)$ be the vertex algebra constructed from 
$A$ and $D$ in Proposition \ref{proposition:comm-alg}
and let $(M,Y_M)$ be a vertex algebra $A$-module.
We call $M$ a {\em vertex algebra $(A,D)$-module} 
to distinguish vertex algebra $A$-modules 
from $\C$-algebra $A$-modules as stated in Section 1.

\begin{proposition}{\rm\cite{B}}\label{proposition:comm-module}
The following hold:
\begin{enumerate}
\item
Let $M$ be a $\C$-algebra $A$-module.
For $a\in A$, define $Y_M(a,x)\in(\End_{\C}M)[[x]]$ by
\begin{align*}
Y(a,x)u&=\sum_{i=0}^{\infty}\dfrac{1}{i!}(D^ia)ux^{i}
\end{align*}
for $u\in M$. Then, $(M,Y_M)$ is a vertex algebra $(A,D)$-module.
\item
Let $(M,Y_M)$ be a vertex algebra $(A,D)$-module such that
$Y(a,x)\in(\End_{\C}M)[[x]]$ for all $a\in A$.
Define an action of $A$ on $M$ by $au=a_{-1}u$ for $a\in A$ and $u\in M$. 
Then, $M$ is a $\C$-algebra $A$-module.
\end{enumerate}
\end{proposition}

By Proposition \ref{proposition:comm-module},
if there exists a vertex algebra $(A,D)$-module $(M,Y_M)$ with $Y_{M}(a,x)\not\in (\End_{\C}M)[[x]]$
for some element $a$ in $A$,
then vertex algebra $(A,D)$-modules and $\C$-algebra $A$-modules are different.
However, no simple criterion for the existence of such a module $(M,Y_M)$ is known.

For a $\C$-linear automorphism $g$ of $V$ of finite order $p$,
set $V^r=\{u\in V\ |\ gu=\zeta_p^{r}u\}, 0\leq r\leq p-1\}$.
We recall the definition of $g$-twisted $V$-modules.
\begin{definition}\label{def:weak-twisted}
A {\em
$g$-twisted $V$-module} $M$ is a vector space equipped with a
linear map
\begin{equation*}
Y_M(\,\cdot\,,x) : V\ni v  \mapsto Y_M(v,x) = \sum_{i \in (1/p)\Z}
v_i x^{-i-1} \in (\End_{\C}M)[[x^{1/p},x^{-1/p}]]
\end{equation*}
which satisfies the following four conditions:

\begin{enumerate}
\item $Y_M(u,x) = \sum_{i \in r/p+\Z}u_i x^{-i-1}$ for $u \in V^r$.
\item $Y_M(u,x)w\in M((x^{1/p}))$ for $u \in V$ and $w \in M$.
\item $Y_M({\mathbf 1},x) = \mathrm{id}_M$.
\item For $u \in V^r$, $v \in V^{s}$, 
$m\in r/T+\Z,\ n\in s/T+\Z$, and $l\in\Z$,
\begin{align*}
&\sum_{i=0}^{\infty}\binom{m}{i}
(u_{l+i}v)_{m+n-i} \\
& = 
\sum_{i=0}^{\infty}\binom{l}{i}(-1)^i
\big(u_{l+m-i}v_{n+i}+(-1)^{l+1}v_{l+n-i}u_{m+i}\big).
\end{align*}
\end{enumerate}
\end{definition}

For a $g$-twisted vertex algebra $(A,D)$-module $(M,Y_{M})$ and 
a linear automorphism $h$ of $A$ which preserves $D$, 
define $(M,Y_{M})\circ h=(M\circ h, Y_{M\circ h})$
by $M\circ h=M$ as vector spaces and $Y_{M\circ h}(a,x) = Y_{M}(ha, x)$
for all $a\in A$.
Then, $(M,Y_{M})\circ h$ is a
$h^{-1}gh$-twisted vertex algebra $(A,D)$-module. 

\section{Finite-dimensional vertex algebra modules over 
fixed point commutative subalgebras}

Throughout this section, $R$ is a connected commutative $\C$-algebra,
$\Sring$ is a commutative $\C$-algebra generated by a single element as an $R$-algebra and
is  a Galois extension of $R$ with Galois group $G$ .
It follows from \cite[Theorem 2.9]{J} that
$A\cong R[Z]/(P(Z))$ as $R$-algebras
for some separable polynomial $P(Z)\in R[Z]$.
Let $D$ be a derivation of $\Sring$ which is invariant under the action of $G$.
For a finite-dimensional vertex algebra $(R,D)$-module $(M,Y_{M})$,
$g\in G$ of order $p$, and 
a linear map $\tilde{Y}(\cdot,x)$ from $\Sring$ to $(\End_{\C}M)((x^{1/p}))$, 
we call $(M,\tilde{Y}_{M})$ a 
{\em
$g$-twisted vertex algebra $(\Sring,D)$-module structure over $(M,Y_{M})$}
if $(M,\tilde{Y}_{M})$ is a  $g$-twisted vertex algebra $(A,D)$-module 
and if $\tilde{Y}(\cdot,x)|_{R}=Y(\cdot,x)$.

In this section, we shall show that 
every finite-dimensional indecomposable vertex algebra $(R,D)$-module 
has a $g$-twisted vertex algebra $(\Sring,D)$-module structure over $(M,Y_M)$
for some $g\in G$.
We use the following notation in \cite[Section 3]{T2}.
For a commutative ring $C$, let $\Mat_n(C)$ denote 
the set of all $n\times n$ matrices with entries in $C$. 
Let $E_n$ denote the $n\times n$ identity matrix
and let $E_{ij}$ denote the matrix whose $(i,j)$ entry is $1$ and all other entries are $0$.
Define $\Delta_{k}(C)=\{(x_{ij})\in\Mat_n(C)\ |\ x_{ij}=0\mbox{ if $i+k\neq j$}\}$ for $0\leq k\leq n$.  
Then, for $a\in \Delta_k(C)$ and $b\in \Delta_{l}(C)$, we have $ab\in \Delta_{k+l}(C)$.
For $X=(x_{ij})\in \Mat_n(C)$ and $k=0,\ldots,n-1$, 
define the matrix $X^{(k)}=\sum_{i=1}^{n}x_{i,i+k}E_{i,i+k}\in \Delta_k(C)$.
For a upper triangular matrix $X$, we see that $X=\sum_{k=0}^{n-1}X^{(k)}$.

Let $A$ be a commutative $\C$-algebra,
$D$ a derivation of $A$, $g$ a $\C$-linear automorphism of $A$ of finite order $p$.
For a vector space $W$ over $\C$ and  a linear map $Y_{W}(\cdot,x)$ from
$A$ to $(\End_{\C}W)[[x^{1/p},x^{-1/p}]]$,
we denote by ${\mathcal A}_{W}(A)$ the subalgebra of $\End_{\C}W$
generated by all coefficients of
$Y_{W}(a,x)$ where $a$ ranges over all elements of $A$.
Let $M$ be a finite-dimensional $g$-twisted vertex algebra $(A,D)$-module. 
Then, ${\mathcal A}_{M}(A)$ is a commutative $\C$-algebra and $M$ is a finite-dimensional ${\mathcal A}_{M}(A)$-module. 
Note that every ${\mathcal A}_{M}(A)$-module becomes $g$-twisted vertex algebra $(A,D)$-module.
Let ${\mathcal J}_{M}(A)$ denote the Jacobson radical of ${\mathcal A}_{M}(A)$.
Recall that the module $\Top M=M/{\mathcal J}_{M}(A)M$ is called {\it the top of }$M$,
which is completely reducible (cf. \cite[Chapter I]{ASS}).
Since ${\mathcal A}_{M}(A)$ is a finite-dimensional commutative $\C$-algebra, 
the Wedderburn--Malcev theorem (cf.\cite[Section 11.6]{P}) says that 
${\mathcal A}_{M}(A)=\oplus_{i=1}^{m}\C e_i\oplus {\mathcal J}_{M}(A)$ 
where $e_1,\ldots,e_m$ are  primitive orthogonal idempotents of ${\mathcal A}_{M}(A)$.
For $U\in {\mathcal A}_{M}(A)((x))$, we denote by $U^{[0]}$
the image of $U$ under the projection 
${\mathcal A}_{M}(A)((x))=\oplus_{i=1}^{m}\C((x))e_i \oplus {\mathcal J}_{M}(A)((x))
\rightarrow\oplus_{i=1}^{m}\C((x))e_i\cong \C((x))^{\oplus m}$.
We denote by $\psi[{A,(M,Y_M)}]$ the $\C$-algebra homomorphism $Y_{M}(\cdot,x)^{[0]}$ 
from $A$ to $\C((x))^{\oplus m}$, which corresponds to the module $\Top M$.
Note that ${\mathcal J}_{M}(A)^{n}((x))=0$, where $n=\dim_{\C}M$.
Since ${\mathcal A}_{M}(A)$ is commutative,
we shall sometimes identify $\End_{\C}M$ with $\Mat_{n}(\C)$
by fixing  a basis of $M$ so that all elements of ${\mathcal A}_{M}(A)$
are upper triangular matrices. Under this identification,
for $U\in {\mathcal A}_{M}(A)((x))$ we see that $U^{[0]}=U^{(0)}$.

Let $M$ be a finite-dimensional indecomposable vertex algebra $(R,D)$-module.
Since ${\mathcal A}_{M}(R)$ is local, we see that
${\mathcal A}_{M}(R)=\C\id \oplus {\mathcal J}_{M}(R)$.
In this case 
we shall often identify the subalgebra $\C((x))\id$ in 
${\mathcal A}_{M}(A)((x))$ with $\C((x))$.
Let $(M,\tilde{Y}_{M})$ be 
a $g$-twisted vertex algebra $(\Sring,D)$-module structure over $(M,Y_M)$.
Since ${\mathcal A}_{M}(R)$ 
is a subalgebra of ${\mathcal A}_{M}(\Sring)$, 
we see that $M$ is an indecomposable ${\mathcal A}_{M}(\Sring)$-module.
Therefore, 
${\mathcal A}_{M}(\Sring)$ is local since 
${\mathcal A}_{M}(\Sring)$ is commutative.
Thus, 
${\mathcal A}_{M}(\Sring)=\C\id\oplus {\mathcal J}_{M}(\Sring)$
and hence $\psi[{\Sring,(M,\tilde{Y}_{M})}]|_{R}=\psi[{R,(M,Y_{M})}]$.
It follows from Nakayama's lemma (cf. \cite[Lemma 2.2]{ASS})
that ${\mathcal J}_{M}(\Sring)M\neq M$ and hence
${\mathcal J}_{M}(\Sring)M={\mathcal J}_{M}(R)M$ is a proper
${\mathcal A}_{M}(\Sring)$-submodule of $M$.
This tells us that $\Top M=M/{\mathcal J}_{M}(R)M$ has 
a $g$-twisted vertex algebra $(\Sring,D)$-module structure over $(\Top M,Y_{\Top M})$.
We conclude that
a $g$-twisted vertex algebra $(\Sring,D)$-module structure 
$(M,\tilde{Y}_M)$ over $(M,Y_M)$ induces
a $g$-twisted vertex algebra $(\Sring,D)$-module structure 
$(\Top M,\tilde{Y}_{\Top M})$ over $(\Top M,Y_{\Top M})$.

Now we state our main theorem.

\begin{theorem}\label{theorem:correspondence}
Let $\Sring$ be a connected commutative $\C$-algebra 
which is a Galois extension of $R$ with Galois group $G$
and let $D$ be a derivation of $\Sring$ which is invariant under the action of $G$.
Suppose $\Sring$ is generated by a single element as an $R$-algebra.
Then, for every non-zero finite-dimensional indecomposable vertex algebra $(R,D)$-module $(M,Y_{M})$, we have the following results:
\begin{enumerate}
\item
$M$ has a $g$-twisted vertex algebra $(\Sring,D)$-module structure over $(M,Y_M)$ for some $g\in G$.
\item 
Let $g\in G$.
If $\Top M$ has a $g$-twisted vertex algebra $(\Sring,D)$-module structure over
$(\Top M, Y_{\Top M})$,
then $M$ has a unique $g$-twisted vertex algebra $(\Sring,D)$-module structure
$(M,\tilde{Y}_{M})$ over $(M,Y_M)$
such that 
$\Top M\cong M/{\mathcal J}_M(\Sring)M$ as $g$-twisted vertex algebra $(\Sring,D)$-modules.
\item Let $g\in G$ and let $(M,\tilde{Y}_{M})$  be a 
$g$-twisted vertex algebra $(\Sring,D)$-module structure over $(M,Y_M)$.
Then, $\tilde{Y}_{M}\circ h, h\in G,$ are all distinct 
homomorphisms from $\Sring$ to $(\End_{\C}M)((x^{1/|g|}))$. 
\item
For each $k=1,2$, let $g_k$ be an element in $G$
and let $(M,\tilde{Y}^{k}_{M})$  be a
$g_k$-twisted vertex algebra $(\Sring,D)$-module structure over $(M,Y_M)$.
Then, $(M,\tilde{Y}^{1}_{M})\circ h\cong (M,\tilde{Y}^{2}_{M})$ for some $h\in G$.
\end{enumerate}
\end{theorem}
\begin{proof}
Set $n=\dim_{\C}M$ and $N=|G|$.
Let the notation be as above.
By \cite[Theorem 2.9]{J},
we may assume $\Sring=R[Z]/(P(Z))$ where
$P(Z)=\sum_{i=0}^{N}P_iZ^i\in R[Z]$ is a separable polynomial.
We denote by $\sR$ the image of the homomorphism $\psi[{R,(M,Y_{M})}] : R\rightarrow \C((x))$,   
by $Q(\sR)$ the quotient field of $\sR$ in $\C((x))$,
by $\theta$ the primitive element $Z+(P(Z))\in \Sring$,
by $\hat{P}(Z)\in({\mathcal A}_{M}(R)((x)))[Z]$ the image of $P(Z)$ under the map $Y_{M}(\cdot,x)$,
and by $\hat{P}^{[0]}(Z)\in\C((x))[Z]$ the image of $P(Z)$ under the map $\psi[{R,(M,Y_{M})}]$.
We write $\hat{P}(Z)=\sum_{i=0}^{N}\hat{P}_i(x)Z^i, \hat{P}_i(x)\in
{\mathcal A}_{M}(R)((x))$.
We use \cite[Lemma 4]{T2} by setting ${\mathcal B}=R\cup\{\theta\}$.

It is well known that any finite extension of $\C((x))$ is $\C((x^{1/j}))$ for some
positive integer $j$ and $\Omega=\cup_{j=1}^{\infty}\C((x^{1/j}))$ is the algebraic closure
of $\C((x))$ (cf. \cite[Corollary 13.15]{E}). 
The field $\C((x^{1/j}))$ becomes a Galois extension of $\C((x))$ whose
Galois group is the cyclic group
generated by the automorphism sending $x^{1/j}$ to $\zeta_jx^{1/j}$.
Let $\sK$ denote 
the splitting field of $\hat{P}^{[0]}(Z)$ in $\Omega$.

\begin{enumerate}
\item[(1)]
Since $\sK$ is a finite extension of $Q(\sR)$ and 
$Q(\sR)$ is a subfield of $\C((x))$,
we see that $\sK\C((x))=\C((x^{1/p}))$ for some positive integer $p$. 
It follows from the isomorphism $\Gal(\C((x^{1/p}))/\C((x)))\cong \Gal(\sK/(\sK\cap\C((x))))$
that $\Gal(\sK/(\sK\cap\C((x))))$ has an element $\sigma$ of order $p$.
Since $\sK$ is a field, there is $\sa\in \sK$ such that $\sigma \sa=\zeta_{p}^{j} \sa$ with $(j,p)=1$.
It follows from $\sa^{p}\in \sK\cap\C((x))$ that
$\sa$ is a root of the polynomial 
$Z^p-\sa^p\in \C((x))[Z]$.
Thus, $\sa$ is an element of $x^{-r/p}\C((x))$ for some integer $r$.
We have $(r,p)=1$ since $\sa^i\not\in \sK^{\langle \sigma\rangle}$ for all $i=1,\ldots,p-1$.
Let $\gamma,\delta$ be integers with $\gamma r+\delta p=1$.
By replacing $\sa$ by $\sa^{\gamma}$,
we have $\sigma \sa=\zeta_{p}^{\gamma j}\sa$ and  $\sa\in x^{-1/p}\C((x))$.
Since $(\gamma j ,p)=1$, by replacing $\sigma$ by a suitable power of $\sigma$,
we have $\sigma \sa=\zeta_{p}\sa$ and  $\sa\in x^{-1/p}\C((x))$.
For all $\bo\in \sK$ with $\sigma \bo=\zeta_{p}^{i}\bo$,
we have $\sigma(\sa^{-i}\bo)=\sa^{-i}\bo$ and hence $\bo\in x^{-i/p}\C((x))$.

Let $T(x)^{[0]}\in \sK$ be a root of $\hat{P}^{[0]}(Z)$.
We have a $\C$-algebra homomorphism $\rho$ from $\Sring=R[Z]/(P(Z))$ to $\sK$ 
with $\rho(\theta)=T(x)^{[0]}$. Since $\sigma $ fixes all elements in $Q(\sR)\subset \sK\cap \C((x))$,
$\sigma (T(x)^{[0]})$ is a root of $\hat{P}^{[0]}(Z)$.
Since $A=R[\theta]$ and
$\rho$ induces a bijection from 
$\{g(\theta)\ |\ g\in G\}$ to the set of all roots of $\hat{P}^{[0]}(Z)$ in $\sK$
as explained just before
Proposition \ref{proposition:comm-alg},
$T(x)^{[0]}$ is a primitive element of $\sK$ over $Q(\sR)$
and there exists a unique $g\in G$ with $\rho(g(\theta))=\sigma(T(x)^{[0]})=\sigma(\rho(\theta))$.
These results tell us that $\rho g=\sigma\rho $ and hence 
the order of $g$ is equal to $p$.

Set $\hat{P}^{[1]}(Z)=\hat{P}(Z)-\hat{P}^{[0]}(Z)\id\in 
{\mathcal J}_{M}(R)((x))[Z]$ and $\hat{P}^{[k]}(Z)=0$ for all $k\geq 2$.
We write
$\hat{P}^{[k]}(Z)=\sum_{i=0}^{N}\hat{P}_{i}(x)^{[k]}Z^i, \hat{P}_{i}(x)^{[k]}\in {\mathcal J}_{M}(R)^{k}((x))$, for all $k\geq 0$.

Since
$\hat{P}^{[0]}(Z)$ has no multiple root in $\Omega$,
we see that $(d\hat{P}^{[0]}/dZ)(T(x)^{[0]})\neq 0$.
For $k=1,2,\ldots,n-1$ we inductively define $T(x)^{[k]}\in {\mathcal J}_{M}(R)^k((x^{1/p}))$ 
by
\begin{align}
T(x)^{[k]}
&=-(\dfrac{d\hat{P}^{[0]}}{dZ}(T(x)^{[0]}))^{-1}\nonumber\\
&\quad\times
\sum_{i=0}^{N}\sum_{j_0=0}^{k}\sum_{
\begin{subarray}{l}
0\leq j_1,\ldots,j_i<k\\
j_0+j_1+\cdots+j_i=k\end{subarray}}
\hat{P}_i(x)^{[j_0]}T(x)^{[j_1]}\cdots T(x)^{[j_i]}\label{eqn:Tk}.
\end{align}
Set $T(x)=\sum_{k=0}^{n-1}T(x)^{[k]}\in{\mathcal A}_{M}(R)((x^{1/p}))$.
Since ${\mathcal J}_{M}(R)^{n}((x))=0$, we have
\begin{align*}
\hat{P}(T(x))&=\hat{P}^{[0]}(T(x)^{[0]})\\
&\quad{}+
\sum_{k=1}^{n-1}\sum_{i=0}^{N}\sum_{
\begin{subarray}{l}
0\leq j_0,j_1,\ldots,j_i\\
j_0+j_1+\cdots+j_i=k\end{subarray}}
\hat{P}_i(x)^{[j_0]}T(x)^{[j_1]}\cdots T(x)^{[j_i]}\\
&=0+
\sum_{k=1}^{n-1}
\Big(T(x)^{[k]}\dfrac{d\hat{P}^{[0]}}{dZ}(T(x)^{[0]})\\
&\quad{}+
\sum_{i=0}^{N}\sum_{j_0=0}^{k}\sum_{
\begin{subarray}{l}
0\leq j_1,\ldots,j_i<k\\
j_0+j_1+\cdots+j_i=k\end{subarray}}
\hat{P}_i(x)^{[j_0]}T(x)^{[j_1]}\cdots T(x)^{[j_i]}\Big)\\
&=0.
\end{align*}
This result enables us to define a homomorphism $\tilde{Y}_{M}(\cdot,x)$ from 
$\Sring= R[Z]/(P(Z))$ to ${\mathcal A}_{M}(R)((x^{1/p}))$ sending $\theta$ to $T(x)$.
Since ${\mathcal A}_{M}(R)$ is commutative,
the subalgebra ${\mathcal A}_{M}(\Sring)$ of $\End_{\C}M$
obtained by $\tilde{Y}_{M}(\cdot,x)$
is commutative. 

For all $b\in \Sring$ with $gb=\zeta_{p}^{i}b$,
we shall show that
$\tilde{Y}_{M}(b,x)\in x^{-i/p}(\End_{\C}M)((x))$.
Set $B(x)=\tilde{Y}_{M}(b,x)$ and $Q(x)=B(x)^{p}\in {\mathcal A}_{M}(R)((x))$.
We identify $\End_{\C}M$ with $\Mat_{n}(\C)$
by fixing  a basis of $M$ so that all elements of ${\mathcal A}_{M}(R)$
are upper triangular matrices. 
We use the expansion $B(x)=\sum_{k=0}^{n-1}B(x)^{(k)},
B(x)^{(k)}\in \Delta_{k}(\End_{\C}M)((x^{1/p}))$.
Since $\zeta_{p}^{i}\rho(b)=\rho(gb)=\sigma (\rho(b))$,
we have already seen that $B(x)^{(0)}=\rho(b)\in x^{-i/p}(\End_{\C}M)((x))$.
By $B(x)^p=Q(x)$, for all $k=1,\ldots,n-1$ we have
\begin{align*}
B(x)^{(k)}
&=-p^{-1}(B(x)^{(0)})^{-p+1}\nonumber\\
&\quad\times (Q(x)^{(k)}+
\sum_{
\begin{subarray}{l}
0\leq j_1,\ldots,j_{p}<k\\
j_1+\cdots+j_{p}=k\end{subarray}}
B(x)^{(j_1)}\cdots B(x)^{(j_p)}).
\end{align*}
It follows by induction on $k$ that
$B(x)^{(k)}\in x^{-i/p}(\End_{\C}M)((x))$ 
and hence $B(x)\in x^{-i/p}(\End_{\C}M)((x))$.

It follows from $P(\theta)=0$
that $0=D(P(\theta))=\sum_{i=0}^{N}(DP_i)\theta^i+(dP/dZ)(\theta)(D\theta)$.
Note that $(d\hat{P}(Z)/dZ)(T(x))$ is an invertible element 
in ${\mathcal A}_{M}(R)((x^{1/p}))$ since $(d\hat{P}^{[0]}(Z)/dZ)(T(x)^{[0]})\neq 0$.
Since $Y_{M}(DP_i,x)=d{Y}_{M}(P_i,x)/dx$ for all $i$,
we have $\tilde{Y}_{M}(D\theta,x)=d\tilde{Y}_{M}(\theta,x)/dx$.

We conclude that 
$(M,\tilde{Y}_{M})$ is a $g$-twisted vertex algebra $(A,D)$-module structure over 
$(M,Y_M)$.

\item[(2)]
We denote the order of $g$ by $p$.
Let $(\Top M,\tilde{Y}_{\Top M})$ be a $g$-twisted vertex algebra $(\Sring,D)$-module structure over $(\Top M,Y_{\Top M})$.
Let us denote by $\varphi$ the map $\tilde{Y}_{\Top M}(\cdot,x) : \Sring\rightarrow \Omega$,
namely $\varphi=\psi[A,(\Top M, \tilde{Y}_{\Top M})]$.
Note that $\varphi|_{R}=\psi[{R,(M,Y_M)}]$ and $\varphi(\theta)$ is a root of $\hat{P}^{[0]}(Z)$ in $\Omega$.
By the same argument as in (1),
we can construct  a root 
$T(x)\in {\mathcal A}_{M}(R)((x^{1/p}))$ of $\hat{P}(Z)$ 
whose semisimple part $T(x)^{[0]}$ is equal to $\varphi(\theta)$. 
The linear homomorphism from $\Sring$ to ${\mathcal A}_{M}(R)((x^{1/p}))$
sending $\theta$ to $T(x)$ induces
a $g$-twisted vertex algebra $(\Sring,D)$-module structure  $(M,\tilde{Y}_{M})$  over 
$(M,Y_M)$.
Since $\theta$ is a primitive element of $\Sring$ over $R$,
we see that $\psi[{\Sring,(M,\tilde{Y}_M)}]=\varphi$.

We shall show the uniqueness of the $g$-twisted vertex algebra $(\Sring,D)$-module structure over 
$(M,Y_M)$ which satisfies the conditions.
Let $(M,\tilde{Y}^1_{M})$ be a 
$g$-twisted vertex algebra $(\Sring,D)$-module structure over 
$(M,Y_M)$ with $\psi[(\Sring,(M,\tilde{Y}^1_M)]=\varphi$.
We identify $\End_{\C}M$ with $\Mat_{n}(\C)$
by fixing  a basis of $M$ so that all elements of ${\mathcal A}_{M}(\Sring)$
are upper triangular matrices. 
Set $U(x)=\tilde{Y}^1_{M}(\theta,x)\in (\Mat_{n}(\C))((x^{1/p}))$. 
We use the expansion $U(x)=\sum_{k=0}^{n-1}U(x)^{(k)}$ and
$\hat{P}_i(x)=\sum_{k=0}^{n-1}\hat{P}_i(x)^{(k)}$,
where 
$U(x)^{(k)},
\hat{P}_i(x)^{(k)}\in \Delta_{k}(\End_{\C}M)((x^{1/p}))$.
Set $\hat{P}^{(0)}(Z)=\sum_{i=0}^{N}\hat{P}_{i}(x)^{(0)}Z^i$.
Under the identification of $\End_{\C}M$ with $\Mat_{n}(\C)$, we have
$\hat{P}^{[0]}(Z)=\hat{P}^{(0)}(Z)$.
Note that $U(x)^{(0)}=\varphi(\theta)$ and we do not assume $U(x)\in {\mathcal A}_{M}(R)((x^{1/p}))$.
We have 
\begin{align*}
0&=\hat{P}(U(x))\\
&=\hat{P}^{(0)}(U(x)^{(0)})\\
&\quad{}+
\sum_{k=1}^{n-1}\sum_{i=0}^{N}\sum_{
\begin{subarray}{l}
0\leq j_0,j_1,\ldots,j_i\\
j_0+j_1+\cdots+j_i=k\end{subarray}}
\hat{P}_i(x)^{(j_0)}U(x)^{(j_1)}\cdots U(x)^{(j_i)}\\
&=0+
U(x)^{(k)}\dfrac{d\hat{P}^{(0)}}{dZ}(U(x)^{(0)})\\
&\quad{}+
\sum_{k=1}^{n-1}\sum_{i=0}^{N}\sum_{j_0=0}^{k}\sum_{
\begin{subarray}{l}
0\leq j_1,\ldots,j_i<k\\
j_0+j_1+\cdots+j_i=k\end{subarray}}
\hat{P}_i(x)^{(j_0)}U(x)^{(j_1)}\cdots U(x)^{(j_i)}
\end{align*}
and hence
\begin{align*}
U(x)^{(k)}
&=-(\dfrac{d\hat{P}^{[0]}}{dZ}(\varphi(\theta))^{-1}\\
&\quad\times
\sum_{i=0}^{N}\sum_{j_0=0}^{k}\sum_{
\begin{subarray}{l}
0\leq j_1,\ldots,j_i<k\\
j_0+j_1+\cdots+j_i=k\end{subarray}}
\hat{P}_i(x)^{(j_0)}U(x)^{(j_1)}\cdots U(x)^{(j_i)}.
\end{align*}
It follows by induction on $k$ that
$U(x)=\sum_{k=0}^{n-1}U(x)^{(k)}$ is uniquely determined by $\varphi(\theta)$ and $\hat{P}(Z)$. 
We conclude that $M$ has a unique $g$-twisted vertex algebra $(\Sring,D)$-module structure
$(M,\tilde{Y}_{M})$ over $(M,Y_M)$
such that $\psi[(\Sring,(M,\tilde{Y}_M)]=\varphi$.

\item[(3)]
Let $h\in G$ with $h\neq 1$.
Since $\theta$ and $h(\theta)$ are distinct roots of $P(Z)$ in $\Sring$,
\cite[Lemma 2.1]{J} says that
$\theta-h(\theta)$ is an invertible element of $\Sring$.
Since $\tilde{Y}_{M\circ h}(\theta,x)=
\tilde{Y}_{M}(h\theta,x)\neq
\tilde{Y}_{M}(\theta,x)$,
we see that $\tilde{Y}_{M\circ h}(\cdot,x)$ is distinct from $\tilde{Y}_{M}(\cdot ,x)$.
This says that 
$\tilde{Y}_{M}\circ h, h\in G,$ are  all distinct 
homomorphisms from $\Sring$ to $(\End_{\C}M)((x^{1/|g|}))$. 
\item[(4)]
For each $k=1,2$, let $g_k$ be an element in $G$ 
and let $(M,\tilde{Y}^{k}_{M})$  be a
$g_k$-twisted vertex algebra $(\Sring,D)$-module structure over $(M,Y_M)$.
We denote $\psi[\Sring,(M,\tilde{Y}^{k}_{M})]$ by $\psi_k$ and 
$\psi[R,(M,Y_M)]$ by $\psi$ briefly.
Since each $\psi_k$ induces a bijection from 
$\{g(\theta)\ |\ g\in G\}$ to the set of all roots of $\hat{P}^{[0]}(Z)$ in $\sK$
as explained just before
Proposition \ref{proposition:comm-alg},
there is an element $h\in G$ with $\psi_1(h(\theta))=\psi_2(\theta)$.
This tells us that
$(\Top M,\tilde{Y}^{1}_{\Top M})\circ h\cong (\Top M,\tilde{Y}^{2}_{\Top M})$
and hence $(M,\tilde{Y}^{1}_{M})\circ h\cong (M,\tilde{Y}^{2}_{M})$ by (2).

\end{enumerate}
\end{proof}

\section{Finite-dimensional vertex algebra $\C[s,s^{-1}]$-modules}

Let 
$\C[s,s^{-1}]$ be the algebra of Laurent polynomials in one variable $s$ over $\C$. 
In this section we shall classify the finite-dimensional vertex algebra $\C[s,s^{-1}]$-modules.
We use the notation introduced in Section 3.
It is easy to see that
every non-zero derivation $D$ of $\C[s, s^{-1}]$
can be expressed as $D=(p(s)/s^{N_q})d/ds$ so that
the polynomials $p(s)$ and  $s^{N_q}$ in $\C[s]$ are coprime.

The following lemma easily follows from \cite[Lemma 4]{T2}.

\begin{lemma}\label{lemma:finite}
Let the notation be as above.
Let $M$ be a finite-dimensional vector space
and let $S(x)=\sum_{i\in \Z}S_{(i)}x^{i}$ be a non-zero element of $(\End_{\C}M)((x))$.
Then, there exists a
vertex algebra $(\C[s,s^{-1}],D)$-module $(M,Y_M)$ 
with $Y_M(s,x)=S(x)$ 
if and only if the following three conditions hold:
\begin{itemize}
\item[(i)] $S(x)$ is an invertible element in $(\End_{\C}M)((x))$.
\item[(ii)] For all $i,j\in \Z$, $S_{(i)}S_{(j)}=S_{(j)}S_{(i)}$.
\item[(iii)] $S(x)^{N_q}dS(x)/dx=p(S(x))$.
\end{itemize}
In this case, for $u(s)\in\C[s,s^{-1}]$ we have $Y_M(u(s),x)=u(S(x))$ and hence
$(M,Y_M)$ is uniquely determined by $S(x)$.
\end{lemma}
\begin{proof}
If
$(M,Y_{M})$ is a vertex algebra $(\C[s,s^{-1}],D)$-module,
then \cite[Lemma 4]{T2} tells us that
the conditions (i)--(iii) are clearly hold and 
$Y_M(u(s),x)=u(S(x))$ for all $u(s)\in\C[s,s^{-1}]$.

Conversely, suppose that $(M,Y_{M})$ satisfies the conditions (i)--(iii).
We use \cite[Lemma 4]{T2} by setting ${\mathcal B}=\{s,s^{-1}\}$.
For $u(s)\in\C[s]$, set $Y_M(u(s),x)=u(S(x))$.
Since $S(x)$ is an invertible element in $(\End_{\C}M)((x))$, this induces a $\C$-algebra homomorphism 
from $\C[s,s^{-1}]$ to $(\End_{\C}M)((x))$.
Since $S(x)^{-1}$ is a polynomial in $S(x)$,
we see that ${\mathcal A}_{M}(\C[s,s^{-1}])$ is commutative.
Since
\begin{align*}
Y_{M}(D(s^{-1}),x)&=
Y_{M}(-(Ds)(s^{-2}),x)\\
&=-Y_{M}(Ds,x)Y_{M}(s,x)^{-2}\\
&=
\dfrac{d}{dx}(Y_{M}(s,x))^{-1}),
\end{align*}
we conclude that $(M,Y_M)$ is a vertex algebra $(\C[s,s^{-1}],D)$-module.
\end{proof}

Let $(M,Y_M)$ be a finite-dimensional indecomposable vertex algebra $\C[s,s^{-1}]$-module.
We identify $\End_{\C}M$ with $\Mat_{n}(\C)$
by fixing  a basis of $M$ so that all elements of ${\mathcal A}_{M}(\C[s,s^{-1}])$
are upper triangular matrices. 
Let $J_n$ denote the following $n\times n$ matrix:
\begin{align*}
J_n&=\begin{pmatrix}
0&1&0&\cdots&0\\
\vdots&\ddots&\ddots&\ddots&\vdots\\
\vdots&&\ddots&\ddots&0\\
\vdots&&&\ddots&1&\\
0&\cdots&\cdots&\cdots&0
\end{pmatrix}.
\end{align*}
We denote $Y_{M}(s,x)$ by $S(x)$.
We use the expansion $S(x)=\sum_{k=0}^{n-1}S(x)^{(k)},
S(x)^{(k)}\in \Delta_{k}(\End_{\C}M)((x))$, as in Section 3.
Recall that under this identification, 
the semisimple part $S(x)^{[0]}$ of $S(x)$ is equal to $S(x)^{(0)}$.  

For all $H(x)=\sum_{i=L}^{\infty}H_ix^i\in(\End_{\C}M)((x))$ with $H_L\neq 0$,
we denote $L$ by $\ld(H(x))$ and $H_{L}$ by $\lc(H(x))$.
Note that if $\ld(S(x)^{[0]})>0$, then 
$\ld((S(x)^{-1})^{[0]})=\ld((S(x)^{[0]})^{-1})<0$.
This implies that if $\ld(S(x)^{[0]})\neq 0$, then
vertex algebra $(\C[s,s^{-1}],D)$-module $(M,Y_M)$ is not a $\C$-algebra $\C[s,s^{-1}]$-module.

\begin{theorem}\label{theorem:untwist}
Let $\alpha$ be a non-zero complex number and $D=(p(s)/s^{N_q})d/ds$ a non-zero derivation of $\C[s, s^{-1}]$
such that the polynomials $p(s)$ and  $s^{N_q}$ of $\C[s]$ are coprime.
We write $p(s)=\sum_{i=L_p}^{N_p}p_is^i$ where $p_{L_p},p_{N_p}$ are non-zero
complex numbers.
Then, the following results hold:
\begin{enumerate}
\item Every finite-dimensional indecomposable vertex algebra $(\C[s,s^{-1}],D)$-module $M$ with $\ld(S(x)^{[0]})=0$
is a $\C$-algebra $A$-module.
\item 
There exists a non-zero finite-dimensional indecomposable  vertex algebra $(\C[s,s^{-1}],D)$-module $M$
with $\ld(S(x)^{[0]})>0$ and with $\lc(S(x)^{[0]})=\alpha$
if and only if $N_q=0$ and $p(0)=\alpha$. Moreover, in this case $\ld(S(x)^{[0]})=1$.
\item There exists a non-zero finite-dimensional indecomposable  vertex algebra $(\C[s,s^{-1}],D)$-module $M$
with $\ld(S(x)^{[0]})<0$ and with $\lc(S(x)^{[0]})=\alpha$
if and only if $N_p=N_q+2$ and $\alpha=-1/p_{N_p}$.
Moreover, in this case $\ld(S(x)^{[0]})=-1$.
\end{enumerate}
In the case of (2) and (3), for every positive integer $n$,
there exists a unique $n$-dimensional indecomposable vertex algebra $(\C[s,s^{-1}],D)$-module
which satisfies the conditions up to isomorphism.
\begin{proof}
We use Lemma \ref{lemma:finite}.
Let $(M,Y_{M})$ be a non-zero finite-dimensional indecomposable vertex algebra $(\C[s,s^{-1}],D)$-module with 
$\lc(S(x)^{[0]})=\alpha$.
Since $M$ is indecomposable, we see that $S(x)^{(0)}\in\C((x))E_n$.
Since $S(x)$ is invertible, we have $S(x)^{(0)}\neq 0$ and
\begin{align}
S(x)^{-1}&
=(S(x)^{(0)}+\sum_{k=1}^{n-1}S(x)^{(k)})^{-1}\nonumber\\
&=\sum_{i=0}^{n-1}(-1)^i(S(x)^{(0)})^{-1-i}(\sum_{k=1}^{n-1}S(x)^{(k)})^i.
\label{eqn:inverse-s}
\end{align}
By Lemma \ref{lemma:finite}, we have
\begin{align}
S(x)^{N_q}\dfrac{dS(x)}{dx}
&=p(S(x)).\label{eqn:expand-s}
\end{align}
and hence
\begin{align}
(S(x)^{(0)})^{N_q}\dfrac{dS(x)^{(0)}}{dx}
&=p(S(x)^{(0)}).\label{eqn:expand-s0}
\end{align}

We shall give a 
formula for $S(x)^{(k)}=\sum_{i\in\Z}S^{(k)}_{(i)}x^i\in \Delta_{k}((\End_{\C}M)((x)))$ for $k=1,2\ldots,n-1$.
By standard Jordan canonical form theory,
we may assume $S_{(0)}=S^{(0)}_{(0)}+S^{(1)}_{(0)}$, that is,
$S^{(j)}_{(0)}=0$ for all $j=2,\ldots,n-1$.
We have the following expansions of $(dS(x)/dx)S(x)^{N_q}$ and $p(S(x))$:
\begin{align*}
&\dfrac{dS(x)}{dx}S(x)^{N_q}\\
&
=\sum_{j_0=0}^{n-1}\dfrac{dS(x)^{(j_0)}}{dx}
\big(
\sum_{
0\leq j_1,\ldots,j_{N_q}\leq n-1}
S(x)^{(j_1)}\cdots S(x)^{(j_{N_q})}\big)\\
&
=
\sum_{
0\leq j_0,j_1,\ldots,j_{N_q}\leq n-1}
\dfrac{dS(x)^{(j_0)}}{dx}
S(x)^{(j_1)}\cdots S(x)^{(j_{N_q})}
\\
&=
\sum_{k=0}^{n-1}
\sum_{\begin{subarray}{c}
0\leq j_0,j_1,\ldots,j_{N_q}\\
j_0+j_1+\cdots+j_{N_q}=k\end{subarray}}
\dfrac{dS(x)^{(j_0)}}{dx}S(x)^{(j_1)}\cdots S(x)^{(j_{N_q})}
\\
&=\dfrac{dS(x)^{(0)}}{dx}(S(x)^{(0)})^{N_q}\\
&\quad{}+
\sum_{k=1}^{n-1}\Big(\dfrac{dS(x)^{(k)}}{dx}(S(x)^{(0)})^{N_q}+
N_q\dfrac{dS(x)^{(0)}}{dx}(S(x)^{(0)})^{N_q-1}S(x)^{(k)}\\
&\quad{}+
\sum_{\begin{subarray}{c}
0\leq j_0,j_1,\ldots,j_{N_q}<k\\
j_0+j_1+\cdots+j_{N_q}=k\end{subarray}}
\dfrac{dS(x)^{(j_0)}}{dx}S(x)^{(j_1)}\cdots S(x)^{(j_{N_q})}\Big)
\end{align*}
and
\begin{align*}
p(S(x))&=
p(S(x)^{(0)})+\sum_{k=1}^{n-1}\big(\dfrac{dp}{ds}(S(x)^{(0)})S(x)^{(k)}\\
&\quad{}+
\sum_{i=0}^{N_p}p_i
\sum_{\begin{subarray}{c}
0\leq j_1,\ldots,j_i<k\\
j_1+\cdots+j_i=k\end{subarray}}
S(x)^{(j_1)}\cdots S(x)^{(j_i)}\big).
\end{align*}
By (\ref{eqn:expand-s}) for $k=1,2,\ldots$, we have a 
formula
\begin{align}
&\dfrac{dS(x)^{(k)}}{dx}\nonumber\\
&=(S(x)^{(0)})^{-N_q}
\Big((-N_q\dfrac{dS(x)^{(0)}}{dx}(S(x)^{(0)})^{N_q-1}
+\dfrac{dp}{ds}(S(x)^{(0)})\big)S(x)^{(k)}\nonumber\\
&\quad{}-
\sum_{\begin{subarray}{c}
0\leq j_0,j_1,\ldots,j_{N_q}<k\\
j_0+\cdots+j_{N_q}=k\end{subarray}}
\dfrac{dS(x)^{(j_0)}}{dx}S(x)^{(j_1)}\cdots S(x)^{(j_{N_q})}\nonumber\\
&\quad{}+
\sum_{i=0}^{N_p}p_i
\sum_{\begin{subarray}{c}
0\leq j_1,\ldots,j_i<k\\
j_1+\cdots+j_i=k\end{subarray}}
S(x)^{(j_1)}\cdots S(x)^{(j_i)}\Big).\label{eqn:semi-nil}
\end{align}
We write $S(x)^{(0)}=\sum_{i=L}^{\infty}S^{(0)}_{(i)}x^i$, where $L=\ld(S(x)^{(0)})$.

Suppose that $L=0$.
We shall show that $S(x)^{(k)}\in(\End_{\C}M)[[x]]$ by induction on $k$.
The case $k=0$ follows from $L=0$.
For $k>0$, suppose that $\ld(S(x)^{(k)})<0$. 
Since $(S(x)^{(0)})^{-N_q}$ is an element of $\C[[x]]$, 
the lowest degree of the right-hand side of (\ref{eqn:semi-nil}) is 
greater than or equal to $\ld(S(x)^{(k)})$ by the induction assumption.
This contradicts that $\ld(dS(x)^{(k)}/dx)=\ld(S(x)^{(k)})-1$. 
It follows from (\ref{eqn:inverse-s}) that $S(x)$ and $S(x)^{-1}$ are elements in $(\End_{\C}M)[[x]]$ and hence
$Y_{M}(a,x)\in (\End_{\C}M)[[x]]$ for all $a\in \C[s,s^{-1}]$.
We conclude that if $L=0$ then $(M,Y_M)$ is a $\C$-algebra
$\C[s,s^{-1}]$-module.
This completes the proof of (1).

Suppose that $L>0$. 
In (\ref{eqn:expand-s0}),
the term with the lowest degree of the left-hand side is
$L(S^{(0)}_{(L)})^{N_q+1}x^{L(N_q+1)-1}$
and the term with the lowest degree of the right-hand side is
$p_{L_p}(S^{(0)}_{(L)})^{L_p}x^{LL_p}$.
Comparing these terms, we have $L(L_p-N_q-1)=-1$
and hence $L=1$ and $L_p=N_q$. We have
$L_p=N_q=0$ since $p(s)$ and $s^{N_q}$ are coprime.
Comparing coefficients of these terms with the lowest degree in (\ref{eqn:expand-s0}), 
we have $D=p(s)d/ds, S^{(0)}_{(1)}=\alpha=p(0)\neq 0$, 
and $S^{(0)}_{(0)}=0$.
For all positive integers $n$,
we shall show the uniqueness of 
$n$-dimensional indecomposable vertex algebra $(\C[s,s^{-1}],D)$-module
which satisfies the conditions in (2).
Setting $N_q=0$ in (\ref{eqn:semi-nil}),
the same argument as in the case of $L=0$ shows that
$S(x)^{(k)}\in(\End_{\C}M)[[x]]$ for
all $k=0,1,\ldots,n-1$.
For all positive integers $m$,
comparing the coefficients of $x^{m}$ in (\ref{eqn:expand-s0}), we have
\begin{align}
(m+1)S^{(0)}_{(m+1)}&=
\sum_{i=0}^{N_p}p_i
\sum_{\begin{subarray}{c}
0\leq j_1,\ldots,j_i\leq m\\
j_0+j_1+\cdots+j_i=m\end{subarray}}
S^{(0)}_{(j_1)}\cdots S^{(0)}_{(j_i)}.
\label{eqn:s0-induc}
\end{align}
It follows by induction on $m$ that every
$S^{(0)}_{(m)}$ is uniquely determined by $S^{(0)}_{(1)}$. 
By (\ref{eqn:semi-nil}) for all $m>0$, $S^{(k)}_{(m)}$ is a polynomial in 
$\{S^{(k)}_{(j)}\ |\ 0\leq j\leq m-1\}
\cup\{S^{(i)}_{(j)}\ |\ 0\leq i\leq k-1, j\geq 0\}$.
Since $S^{(i)}_{(0)}=0$ for all $i=2,\ldots,n-1$,
it follows by induction on $k$ and $m$ 
that every $S^{(k)}_{(m)}$ is a polynomial in $S^{(1)}_{(0)}$
and hence is uniquely determined by $S^{(1)}_{(0)}$.
Since $S^{(1)}_{(0)}$ is the nilpotent part of $S_{(0)}$ and 
$M$ is indecomposable, $S^{(1)}_{(0)}$ conjugates to $J_n$.
Thus, we have shown that the uniqueness of 
$n$-dimensional indecomposable vertex algebra $(\C[s,s^{-1}],D)$-module
which satisfies the conditions in (2).

Conversely, suppose that $\alpha=p(0)$.
Set $S^{(0)}_{(1)}=\alpha$ and $S^{(0)}_{(i)}=0$ for all non-positive integers $i$.
By (\ref{eqn:s0-induc}) we can inductively define $S^{(0)}_{(m)}$ for $m=2,3,\ldots$. 
The obtained $S(x)^{(0)}=\sum_{i=1}^{\infty}S_{(i)}^{(0)}x^i\in \C[[x]]$
satisfies $\ld(S(x)^{(0)})=1$, $\lc(S(x)^{(0)})=\alpha$, and (\ref{eqn:expand-s0}).
Set $S^{(1)}_{(0)}=J_n$,
$S^{(k)}_{(0)}=0$ for all $k=2,\ldots,n-1$,
and $S^{(k)}_{(i)}=0$ for all $k=1,\ldots,n-1$ and all negative integers $i$.
After (\ref{eqn:s0-induc}),
we have seen that every $S^{(k)}_{(m)}$ is a polynomial in $S^{(1)}_{(0)}$ if it exists.
By the same argument, we can inductively define $S^{(k)}_{(m)}\in \End_{\C}M$
for $k=1,2,\ldots,n-1$ and $m=1,2,\ldots$.
By the argument to get (\ref{eqn:semi-nil}) and (\ref{eqn:s0-induc}) above,
it is easy to see that the obtained $S(x)=\sum_{k=0}^{n-1}S(x)^{(k)}\in (\End_{\C}M)[[x]]$ 
satisfies (\ref{eqn:expand-s}).
Since all coefficients of $S(x)$ are polynomials 
in $S^{(1)}_{(0)}=J_n$, we see that 
$S_{(i)}S_{(j)}=S_{(j)}S_{(i)}$ for all $i,j\in\Z$.
Thus, we have an $n$-dimensional vertex algebra $(\C[s,s^{-1}],D)$-module $M$
with $\ld(S(x)^{(0)})=1$ and with $\lc(S(x)^{(0)})=\alpha$.
This completes the proof of (2).

Next suppose that $L<0$.
Set $\tilde{s}=1/s$. By (\ref{eqn:inverse-s}), we have $Y_{M}(\tilde{s},x)^{[0]}=(S(x)^{-1})^{[0]}=(S(x)^{[0]})^{-1}$ and 
\begin{align*}
D&=-\tilde{s}^{N_q+2}p(1/\tilde{s})\dfrac{d}{d\tilde{s}}.
\end{align*}
Since $S(x)^{-1}$ is a polynomial in $S(x)$,
all coefficients in $S(x)^{-1}$ are commutative.
Thus, this case reduces to the case of $L>0$.
\end{proof}
\end{theorem}

\section{Examples}
Throughout this section,
$D$ is a non-zero derivation of $\C[s,s^{-1}]$. 
For a positive integer $n$,
the $\C$-algebra $A=\C[s,s^{-1}][t]/(t^n-s)$ is a 
Galois extension of $\C[s,s^{-1}]$
(cf. \cite[Lemma 5.1 in Chapter 0]{G}).
The Galois group of $A$ over $\C[s,s^{-1}]$
is the cyclic group of order $n$ generated by
$\tau$ with $\tau(t)=\zeta_n t$.
Since $t^n-s$ is an irreducible element
in the unique factorization domain $\C[s,s^{-1}][t]$,
$t^n-s$ is a prime element.
Therefore, $A$ is an integral domain and hence is a connected $\C$-algebra.
We can extend $D$ to a unique derivation of $A$, which we denote 
by the same notation $D$, by setting
$D(t)=s^{-1}tD(s)/n$.
It is easy to see that $D$ is invariant under the action of $\tau$.

In Theorem \ref{theorem:untwist},
we have classified the finite-dimensional indecomposable $(\C[s,s^{-1}],D)$-modules
$(M,Y_M)$ which are not $\C$-algebra $\C[s,s^{-1}]$-modules.
In this section, we shall investigate twisted 
vertex algebra $(A,D)$-module structures over $(M,Y_M)$.
We denote $Y_{M}(s,x)$ by $S(x)$ and 
$S(x)^{[0]}=\sum_{i=L}^{\infty}S^{[0]}_{(i)}x^i\in\C((x))$ with
$S^{[0]}_{(L)}\neq 0$ as in Section 4.
It follows from Theorem \ref{theorem:untwist} that $L=\ld(S(x)^{[0]})=1$ or $-1$. 
\begin{proposition}
Let $(M,Y_M)$ be a finite-dimensional indecomposable 
vertex algebra $(\C[s,s^{-1}],D)$-module which 
is not a $\C$-algebra $A$-module.
Set $L=\ld(S(x)^{[0]})$.
Then, for the $\C$-algebra $A=\C[s,s^{-1}][t]/(t^n-s)$,
$(M,Y_M)$ has exactly $n$ $\tau^{-L}$-twisted vertex algebra $(A,D)$-module structure.
\end{proposition}
\begin{proof}
We use the notation in the proof of Theorem \ref{theorem:correspondence} (1).
If $\ld(S(x)^{[0]})=1$,
then every root of the polynomial $Z^n-S(x)^{[0]}$ in $\Omega=\cup_{i=1}^{\infty}\C((x^{1/i}))$
is an element in $x^{1/n}\C((x))=x^{-(-1/n)}\C((x))$.
It follows from the argument in the proof of Theorem \ref{theorem:correspondence} (1) that 
$(M,Y_M)$ has a $\tau^{-1}$-twisted vertex algebra $(A,D)$-module structure
$(M,\tilde{Y}_M)$ with $\tilde{Y}_{M}(t,x)\in x^{-(-1/n)}(\Mat_{\C} M)((x))$.
We conclude by Theorem \ref{theorem:correspondence} (3) and (4) that 
$(M,Y_M)$ has exactly $n$ $\tau^{-1}$-twisted vertex algebra $(A,D)$-module structures.
The same argument tells us that
if $\ld(S(x)^{[0]})=-1$, then
$(M,Y_M)$ has exactly $n$ $\tau$-twisted vertex algebra $(A,D)$-module structures.
\end{proof}

\end{document}